\documentclass[11pt]{article}
\usepackage{graphicx}
\usepackage{color}
\usepackage{amsmath}
\usepackage{amssymb}
\usepackage{amscd}
\usepackage{framed}
\usepackage{bbm}

\usepackage{fancyhdr,a4wide}
\usepackage{amsthm}

\newcommand{\R}{\mathbb{R}}
\newcommand{\inr}[1]{\left\langle #1 \right\rangle}

\newcommand{\E}{\mathbb{E}}

\newcommand{\eps}{\varepsilon}

\newtheorem{Theorem}{Theorem}[section]
\newtheorem{Lemma}[Theorem]{Lemma}

\newtheorem{Definition}[Theorem]{Definition}

\newtheorem{Corollary}[Theorem]{Corollary}
\newtheorem{Remark}[Theorem]{Remark}

\newtheorem{Question}[Theorem]{Question}

\numberwithin{equation}{section}

\def \proof {\noindent {\bf Proof.}\ \ }

\def \endproof
{{\mbox{}\nolinebreak\hfill\rule{2mm}{2mm}\par\medbreak}}
\def\IND{\mathbbm{1}}

\def\IND{\mathbbm{1}}

\begin{document}
\title{On the geometry of random polytopes}
\author{
Shahar Mendelson \thanks{LPSM, Sorbonne University,  and Mathematical Sciences Institute, The Australian National University. Email: shahar.mendelson@upmc.fr}}

\maketitle

\begin{abstract}
We present a simple proof to a fact recently established in \cite{GLT}: let $\xi$ be a symmetric random variable that has variance $1$, let $\Gamma=(\xi_{ij})$ be an $N \times n$ random matrix whose entries are independent copies of $\xi$, and set $X_1,...,X_N$ to be the rows of $\Gamma$. Then under minimal assumptions on $\xi$ and as long as $N \geq c_1n$, 
$$
c_2 \bigl(B_\infty^n \cap \sqrt{\log(eN/n)} B_2^n \bigr) \subset {\rm absconv}(X_1,...,X_N)
$$
with high probability. 
\end{abstract}

\section{Introduction}
Let $\xi$ be a symmetric random variable that has variance $1$ and let $X=(\xi_1,...,\xi_n)$ be the random vector whose coordinates are independent copies of $\xi$. Consider a random matrix $\Gamma$ whose rows $X_1,...,X_N$  are independent copies of $X$. In this note we explore the geometry of the random polytope
$$
K={\rm absconv}(X_1,...,X_N) = \Gamma^* B_1^N;
$$
specifically, we study whether $K$ is likely to contain a large canonical convex body.
\vskip0.4cm
One of the first results in this direction is from \cite{G}, where it is shown that if $\xi$ is the standard gaussian random variable, $0<\alpha<1$ and $N \geq c_0(\alpha) n$, then
\begin{equation} \label{eq:gaussian-est}
c_1(\alpha) \sqrt{\log(eN/n)} B_2^n \subset {\rm absconv}(X_1,...,X_N)
\end{equation}
with probability at least $1-2\exp(-c_2N^{1-\alpha} n^{\alpha})$. It should be noted that this estimate cannot be improved---up to the dependence of the constants on $\alpha$ (see, for example, the discussion in Section~4 of \cite{LPRT}).

\vskip0.4cm

The proof of \eqref{eq:gaussian-est} relies heavily on the tail behaviour of the gaussian random variable. It is therefore natural to try and extend \eqref{eq:gaussian-est} beyond the gaussian case, to random polytopes generated by more general random variables that still have `well-behaved' tails. The optimal \emph{subgaussian estimate} was established in \cite{LPRT}:
\begin{Theorem} \label{thm:LPRT} 
Let $\xi$ be a mean-zero random variable that has variance $1$ and is $L$-subgaussian\footnote{A centred random variable is $L$-subgaussian if for every $p \geq 2$, $\|\xi\|_{L_p} \leq L \sqrt{p}\|\xi\|_{L_2}$.}. Let $0<\alpha<1$ and set $N \geq c_0(\alpha)n$. Then with probability at least $1-2\exp(-c_1N^{1-\alpha}n^{\alpha})$
\begin{equation} \label{eq:subgaussian-est}
c_2(\alpha)\bigl(B_\infty^n \cap \sqrt{\log(eN/n)} B_2^n\bigr) \subset {\rm absconv}(X_1,...,X_N),
\end{equation}
where $c_0$ and $c_2$ are constants that depend on $\alpha$ and $c_1$ is an absolute constant.
\end{Theorem}

\begin{Remark}
Note that the body ${\rm absconv}(X_1,...,X_N)$ contains in \eqref{eq:subgaussian-est} is slightly smaller than in \eqref{eq:gaussian-est}, as one has to intersect the Euclidean ball from \eqref{eq:gaussian-est} with the unit cube.
\end{Remark}

While Theorem \ref{thm:LPRT} resolves the problem when $\xi$ is subgaussian, the situation is less clear when $\xi$ is heavy-tailed. That naturally leads to the following question:

\begin{Question} \label{Qu:main}
Under what conditions on $\xi$ one still has that for $N \geq c_1 n$,
\begin{equation} \label{eq:in-Qu:main}
c_2(B_\infty^n \cap \sqrt{\log(eN/n)} B_2^n) \subset {\rm absconv}(X_1,...,X_N)
\end{equation}
with high probability?
\end{Question}

Following the progress in \cite{KKR}, where Question \ref{Qu:main} had been studied under milder moment assumptions on $\xi$ than in Theorem \ref{thm:LPRT}, Question \ref{Qu:main} was answered in \cite{GLT} under a minimal small-ball condition on $\xi$.
\begin{Definition} \label{def:small-ball}
A mean-zero random variable $\xi$ satisfies a small-ball condition with constants $\kappa$ and $\delta$ if
\begin{equation} \label{eq:SB}
Pr(|\xi| \geq \kappa) \geq \delta.
\end{equation}
\end{Definition}

\begin{Theorem} \label{thm:main} \cite{GLT}
Let $\xi$ be a symmetric, variance $1$ random variable that satisfies \eqref{eq:SB} with constants $\kappa$ and $\delta$. For $0<\alpha<1$ there are constants $c_1,c_2$ and $c_3$ that depend on $\kappa, \delta$ and $\alpha$ for which the following holds. If $N \geq c_1 n$ then with probability at least $1-2\exp(-c_2N^{1-\alpha}n^{\alpha})$,
$$
c_3 (B_\infty^n \cap \sqrt{\log(eN/n)} B_2^n) \subset {\rm absconv}(X_1,...,X_N).
$$
\end{Theorem}

\begin{Remark}
The assumption made in \cite{GLT} is slightly stronger than in Theorem \ref{thm:main}; namely, that for every $x \in \R$, $Pr(|\xi-x| \geq \kappa) \geq \delta$. However, \eqref{eq:SB} suffices for the proof. At the same time, in \cite{GLT} the random variables $(\xi_{ij})$ are only assumed to be independent, symmetric and variance $1$, with each one of the $\xi_{ij}$'s satisfying \eqref{eq:SB} with the same constants $\kappa$ and $\delta$. In what follows we consider only the case in which $\xi_{ij}$ are independent copies of a single random variable $\xi$---though extending the presentation to the independent case  is straightforward.
\end{Remark}

The original proof of Theorem \ref{thm:main} is based on the construction of a well-chosen net, and that construction is rather involved. Here we present a much simpler argument that is based on the \emph{small-ball method} (see, e.g., \cite{Men-GAFALect, Men-ACM,Men-GL}). As an added value, the method presented here gives more information than the assertion of Theorem \ref{thm:main}, as is explained in what follows.

\vskip0.4cm

The starting point of the proof of Theorem \ref{thm:main} is straightforward: let
$$
K={\rm absconv}(X_1,...,X_n) = \Gamma^* B_1^N
$$
and set
$$
L= \bigl(B_\infty^n \cap \sqrt{\log(eN/n)} B_2^n\bigr).
$$
By comparing the support functions of $L$ and of $K$, one has to show that with the wanted probability, for every $z \in \R^n$, $h_L(z) \leq h_{cK}(z)$. And, since $h_{cK}(z)=c\|\Gamma z\|_\infty$, Theorem \ref{thm:main} can be established by showing that for suitable constants $c_0$ and $c_1$,
\begin{equation} \label{eq:what-we-show}
Pr( \exists z \in \partial L^\circ \ \ \ \|\Gamma z\|_\infty \leq c_0 ) \leq 2\exp(-c_1 N^{1-\alpha}n^{\alpha}).
\end{equation}

What we actually show is a stronger statement than \eqref{eq:what-we-show}: not only is there a high probability event on which 
$$
\inf_{z \in \partial L^\circ} \|\Gamma z\|_\infty \geq c_0,
$$
but in fact, on that ``good event", for each $z \in \partial L^\circ$, $\Gamma z$ has $\sim N^{1-\alpha}n^{\alpha}$ large coordinates, with each one of these coordinates satisfying that $|\inr{z,X_i}| \geq c_0$. Thus, the fact that $\|\Gamma z\|_\infty \geq c_0$ is exhibited by many coordinates and not just by a single one.

\vskip0.3cm

Proving that indeed, with high probability the smallest cardinality
$$
\inf_{z \in \partial L^\circ} |\{i : |\inr{z,X_i}| \geq c_0\}|
$$
is large is carried out in two steps:

\vskip0.4cm
\noindent{\bf Controlling a single point.} For $0<\alpha<1$ and a well chosen $c_0=c_0(\alpha)$ one establishes an \emph{individual estimate}: that for every fixed $z \in \partial L^\circ$,
$$
Pr(|\inr{z,X}| \geq 2c_0) \geq 4\left(\frac{n}{N}\right)^{\alpha}.
$$
In particular, if $X_1,...,X_N$ are independent copies of $X$ then with probability at least $1-2\exp(-c_2 N^{1-\alpha} n^{\alpha})$,
\begin{equation} \label{eq:card-intro}
\bigl|\{i : |\inr{z,X_i}| \geq 2c_0\} \bigr| \geq 2N^{1-\alpha} n^{\alpha}.
\end{equation}
\vskip0.4cm
\noindent{\bf From a single function to uniform control.} Thanks to the high probability estimate with which \eqref{eq:card-intro} holds, it is possible to control uniformly any subset of $\partial L^\circ$ whose cardinality is at most $\exp(c_2 N^{1-\alpha} n^{\alpha}/2)$. Let ${\cal T}$ be a minimal $\rho$-cover of $\partial L^\circ$ with respect to the $\ell_2$ norm of the allowed cardinality. For every $z \in \partial L^\circ$, let $\pi z \in {\cal T}$ that satisfies $\|z-\pi z\|_2 \leq \rho$. The wanted uniform control is achieved by showing that
$$
\sup_{z \in \partial L^\circ} \bigl|\{i : |\inr{z-\pi z, X_i}| \geq c_0\} \bigr| \leq N^{1-\alpha} n^{\alpha}
$$
with probability at least $1-2\exp(-c_3(\alpha) N^{1-\alpha} n^{\alpha})$.

Indeed, combining the two estimates it follows that with probability at least
$$
1-2\exp(-c(\alpha)N^{1-\alpha}n^{\alpha}),
$$
for every $z \in \partial L^\circ$, one has that
$$
\bigl|\{ i : |\inr{\pi z,X_i}| \geq 2c_0\}| \geq 2N^{1-\alpha} n^{\alpha}
$$
and
$$
\bigl|\{ i : |\inr{z-\pi z,X_i}| \geq c_0 \}\bigr| \leq N^{1-\alpha}n^{\alpha}.
$$
Hence, on that event, for every $z \in \partial L^\circ$ there is $J_z \subset \{1,...,n\}$ of cardinality at least $N^{1-\alpha} n^{\alpha}$, and for every $j \in J_z$,
$$
|\inr{z,X_i}| \geq |\inr{\pi z,X_i}| - |\inr{z-\pi z,X_i}| \geq c_0,
$$
implying that
$$
\inf_{z \in \partial L^\circ} \bigl|\{i: |\inr{z,X_i}| \geq c_0 \} \bigr| \geq N^{1-\alpha} n^{\alpha};
$$
in particular, $\inf_{z \in \partial L^\circ} \|\Gamma z\|_{\infty} \geq c_0$ as required.

\vskip0.4cm
In the next section this line of reasoning is used to prove Theorem \ref{thm:main}.

\section{Proof of Theorem \ref{thm:main}}
Before we begin the proof, let us introduce some notation. Throughout, absolute constant are denoted by $c,c_1,c^\prime$ etc. . Unless specified otherwise, the value of these constants may change from line to line. Constants that depend on some parameter $\alpha$ are denoted by $c(\alpha)$. We write $a \lesssim b$ if there is an absolute constant $c$ such that $a \leq cb$; $a \lesssim_\alpha b$ implies that $a \leq c(\alpha)b$; and $a \sim b$ if both $a \lesssim b$ and $b \lesssim a$.

\vskip0.4cm

The required estimate for a single point follows very closely  ideas from \cite{MS}, which had been developed for obtaining lower estimates on the tails of marginals of the Rademacher vector $(\eps_i)_{i=1}^n$, that is, on
$$
Pr(|\sum_{i=1}^n \eps_i z_i| > t)
$$
as a function of the `location' in $\R^n$ of $(z_i)_{i=1}^n$.

Fix $1 \leq r \leq n$ and consider the interpolation body $L_r=B_\infty^n \cap \sqrt{r} B_2^n$ and its dual $L_r^\circ = {\rm conv}(B_1^n \cup (1/\sqrt{r})B_2^n)$. The key estimate one needs to establish the wanted individual control is:
\begin{Theorem} \label{thm:single}
There exist constants $c^\prime$ and $c^{\prime \prime}$ that depend only on the small-ball constants of $\xi$ ($\kappa$ and $\delta$) such that if $z \in \partial L_r^\circ$ then $$
Pr \bigl( |\inr{z,X}| \geq c^\prime \bigr) \geq 2\exp(-c^{\prime \prime} r).
$$
\end{Theorem}

The proof of Theorem \ref{thm:single} is based on some well-known facts on the interpolation norm $\| \ \|_{L_r^\circ}$.
\begin{Lemma} \label{lemma:trivial-1}
There exists an absolute constant $c_0$ such that for every $z \in \R^n$,
$$
\|z\|_{L_r^\circ} \leq \sum_{i=1}^r z_i^* + \sqrt{r}\bigl(\sum_{i >r } (z_i^2)^*\bigr)^{1/2} \leq c_0 \|z\|_{L_r^\circ},
$$
where $(z_i^*)_{i=1}^n$ is the nonincreasing rearrangement of $(|z_i|)_{i=1}^n$.

Moreover, for very $z \in \R^n$ there is a partition of $\{1,...,n\}$ to $r$ disjoint blocks $I_1,...,I_r$ such that
$$
\frac{\|z\|_{L_r^\circ}}{\sqrt{2}} \leq \sum_{j=1}^r \bigl(\sum_{i \in I_j} z_i^2\bigr)^{1/2} \leq \|z\|_{L_r^\circ}.
$$
\end{Lemma}

The first part of Lemma \ref{lemma:trivial-1} is due to Holmstedt (see Theorem~4.1 in \cite{H}) and it gives useful intuition on the nature of the norm $\| \ \|_{L_r^\circ}$. The second part is Lemma~2 from \cite{MS} and it plays an essential role in what follows.

\vskip0.4cm

Before proving Theorem \ref{thm:single}, we require an additional observation that is based on the small-ball condition satisfies by $\xi$.
\begin{Lemma} \label{lemma:L-1-est}
Let $J \subset \{1,...,n\}$ and set $Y=\sum_{j \in J} z_j \xi_j$. Then
$$
\E|Y| \geq c(\kappa,\delta) \bigl(\sum_{j \in J} z_j^2\bigr)^{1/2},
$$
where $c(\kappa,\delta)<1$ is a constant the depends only on $\xi$'s small-ball constants $\kappa$ and $\delta$.
\end{Lemma}

\proof
Let $(\eps_j)_{j \in J}$ be independent, symmetric, $\{-1,1\}$-valued random variables that are also independent of $(\xi_j)_{j \in J}$. Recall that $\xi$ is symmetric and therefore $(\xi_j)_{j \in J}$ has the same distribution as $(\eps_j \xi_j)_{j \in J}$. By Khintchine's inequality it is straightforward to verify that
$$
\E |Y| = \E_\xi \E_\eps \bigl|\sum_{j \in J} \eps_j z_j \xi_j\bigr| \gtrsim \E_\xi \bigl(\sum_{j \in J} z_j^2 \xi_j^2\bigr)^{1/2}.
$$
Let $(\eta_j)_{j \in J} = \IND_{ \{|\xi_j| \geq \kappa\}}$; thus, the $\eta_j$'s are iid $\{0,1\}$-valued random variables whose mean is at least $\delta$, and point-wise
$$
\bigl(\sum_{j \in J} z_j^2 \xi_j^2\bigr)^{1/2} \geq \kappa \bigl(\sum_{j \in J} \eta_j z_j^2 \bigr)^{1/2}.
$$
Hence, and all that is left to complete the proof is to show that 
$$
\E \bigl(\sum_{j \in J} \eta_j z_j^2 \bigr)^{1/2} \geq c(\delta) \bigl(\sum_{j \in J} z_j^2 \bigr)^{1/2}.
$$

Let $a_j=z_j^2/(\sum_{j \in J} z_j^2)$ and in particular, $\|(a_j)_{j \in J}\|_1 =1$. Assume without loss of generality that $J=\{1,...,\ell\}$ and that the $a_j$'s are non-increasing, let $\gamma>0$ be a parameter to be specified in what follows, and set $p=\E \eta_1 \geq \delta$. 

Consider two cases:

\noindent{$\bullet$} If $a_1 \geq \gamma p$ then with probability at least $p$, $\sum_{j=1}^\ell \eta_j a_j \geq a_1 \geq \gamma p$. In that case
    $$
    \E \bigl(\sum_{j=1}^\ell \eta_j a_j\bigr)^{1/2} \geq \sqrt{\gamma} p^{3/2} \geq \sqrt{\gamma} \delta^{3/2}.
    $$

\noindent{$\bullet$} Alternatively, $a_1 \leq \gamma p$, implying that
$$
A=\sum_{j=1}^\ell a_j^2 \leq a_1 \sum_{j=1}^\ell a_j \leq \gamma p
$$
because $\|(a_j)_{j=1}^\ell\|_1 =1$. 

By Bernstein's inequality,
\begin{equation*}
Pr \Bigl( \bigl|\sum_{j=1}^\ell (\eta_j-p)a_j\bigl| \geq \frac{p}{2} \Bigr) \leq 2\exp\Bigl(-c_0 \min \Bigl\{ \frac{(p/2)^2}{pA},\frac{p/2}{a_1}\Bigr\} \Bigr)
\leq 2 \exp(-c_1/\gamma) \leq \frac{1}{2}
\end{equation*}
provided that $\gamma$ is a small-enough absolute constant. Using, once again, that $\|(a_j)_{j=1}^\ell\|_1 =1$ it is evident that  with probability $1/2$, $\sum_{j=1}^\ell \eta_j a_j \geq (1/2) p$ and therefore
$$
\E \bigl(\sum_{j=1}^\ell \eta_j a_j \bigr)^{1/2} \geq \frac{\sqrt{p}}{4} \geq \frac{\sqrt{\delta}}{4}.
$$

Thus, setting $c(\kappa,\delta) \sim \kappa \delta^{3/2}$ one has that 
$$
\bigl(\sum_{j=1}^\ell  z_j^2 \xi_j^2\bigr)^{1/2} \geq c(\kappa,\delta) \bigl(\sum_{j=1}^\ell z_j^2 \bigr)^{1/2},
$$
as claimed.

\endproof

\noindent{\bf Proof of Theorem \ref{thm:single}.} Fix $z \in \partial L_r^\circ$ and recall that by Lemma \ref{lemma:trivial-1}
there is a decomposition of $\{1,...,n\}$ to disjoint blocks $(I_j)_{j=1}^r$ such that
\begin{equation} \label{eq:sum-of-l2}
\sum_{j=1}^r \bigl(\sum_{i \in I_j} z_i^2 \bigr)^{1/2} \geq \frac{1}{\sqrt{2}}.
\end{equation}
Let $Y_j = \sum_{i \in I_j} z_i \xi_i$; observe that $Y_1,...,Y_r$ are independent random variables and that by Lemma \ref{lemma:L-1-est},
$$
\E |Y_j| \geq c(\kappa,\delta) \bigl(\sum_{i \in I_j} z_i^2 \bigr)^{1/2}
$$
for a constant $0<c(\kappa,\delta)<1$.

At the same time,
$$
\E |Y_j|^2  = \sum_{i \in I_j} z_i^2 \E \xi_i^2 = \sum_{i \in I_j} z_i^2.
$$
Therefore, by the Paley-Zygmund inequality (see, e.g., \cite{dPG}), for any $0<\theta<1$,
$$
Pr(|Y_j| \geq \theta \E|Y_j|) \geq (1-\theta^2)\frac{(\E |Y_j|)^2}{\E Y_j^2}.
$$
Setting $\theta=1/2$,
$$
Pr\Bigl(|Y_j| \geq  \frac{1}{2} c(\kappa,\delta)  \bigl(\sum_{i \in I_j} z_i^2 \bigr)^{1/2} \Bigr) \geq \frac{3}{4} c^2(\kappa,\delta),
$$
and since $Y_j$ is a symmetric random variable (because the $\xi_i$'s are symmetric), it follows that
$$
Pr\Bigl(Y_j \geq  \frac{1}{2} c(\kappa,\delta)  \bigl(\sum_{i \in I_j} z_i^2 \bigr)^{1/2} \Bigr) \geq \frac{3}{8} c^2(\kappa,\delta) \equiv c_1(\kappa,\delta).
$$
For $1 \leq j \leq r$ let
$$
{\cal B}_j=\Bigl\{Y_j \geq  \frac{1}{2} c(\kappa,\delta)  \bigl(\sum_{i \in I_j} z_i^2 \bigr)^{1/2}\Bigr\}
$$
which are independent events. Hence,
\begin{align*}
Pr \Bigl( \sum_{i=1}^n \xi_i z_i  \geq  \frac{1}{2} c(\kappa,\delta) \sum_{j=1}^r \bigl(\sum_{i \in I_j} z_i^2 \bigr)^{1/2} \Bigr) = & Pr\Bigl( \sum_{j=1}^r Y_j \geq \frac{1}{2} c(\kappa,\delta) \sum_{j=1}^r \bigl(\sum_{i \in I_j} z_i^2 \bigr)^{1/2} \Bigr)
\\
\geq & \prod_{j=1}^r  Pr  ({\cal B}_j) \geq c_1^r(\kappa,\delta).
\end{align*}
Thus, by \eqref{eq:sum-of-l2}, if $c^\prime = \frac{1}{4} c(\kappa,\delta)$ and $c^{\prime \prime}=\log(1/c_1(\kappa,\delta))>0$, one has
$$
Pr\bigl( \sum_{i=1}^n \xi_i z_i \geq c^\prime\bigr) \geq \exp(-c^{\prime \prime}r).
$$
\endproof

From here on, the constants $c^\prime$ and $c^{\prime \prime}$ denote the constants from Theorem \ref{thm:single}.

\begin{Corollary} \label{cor:single}
For $0<\alpha < 1$, $\kappa$ and $\delta$ there are constants $c_0$ and $c_1$ that depend on $\alpha$, $\kappa$ and $\delta$, and an absolute constant $c_2$ for which the following holds. If $N \geq c_0n$, $r \leq c_1 \sqrt{\log(eN/n)}$ and $z \in \partial L_r^\circ$ then with probability at least $1-2\exp(-c_2 N^{1-\alpha}n^{\alpha})$,
$$
\bigl| \{i : |\inr{z,X_i}| \geq c^\prime \} \bigr| \geq 2N^{1-\alpha} n^{\alpha}.
$$
\end{Corollary}

\proof
Let $z \in \partial L_r^\circ$, and invoking Theorem \ref{thm:single},
$$
Pr \bigl( |\inr{z,X} | \geq c^\prime \bigr) \geq \exp(-c^{\prime \prime} r)
$$
where $c^\prime$ and $c^{\prime \prime}$ depend only on $\kappa$ and $\delta$.

Let $r_0=c_1 \log(eN/n)$ such that $\exp(-c^{\prime \prime} r_0) \geq 4(n/N)^{\alpha}$; thus, $c_1=c_1(\alpha,\kappa,\delta)$. If $r \leq r_0$, $X_1,...,X_N$ are independent copied of $X$ and $\eta_i = \IND_{\{ |\inr{z,X_i}| \geq c^\prime\}}$,  then $\E \eta_i \geq 4(n/N)^{\alpha}$. Hence, by a standard concentration argument, with probability at least $1-2\exp(-c_2 N^{1-\alpha} n^{\alpha})$,
$$
\bigl| \{i : |\inr{z,X_i}| \geq c^\prime\} \bigr| \geq 2N^{1-\alpha}n^{\alpha},
$$
where $c_2$ is an absolute constant.
\endproof

Thanks to the high probability estimate with which Corollary \ref{cor:single} holds, one can control uniformly all the elements of a set ${\cal T} \subset  \partial L_r^\circ$ as long as $|{\cal T}| \leq \exp(c_0 N^{1-\alpha}n^{\alpha})$ for a suitable absolute constant $c_0$, and as long as $r \leq c(\alpha,\kappa,\delta) \log(eN/n)$. In that case, there is an event of probability at least $1-2\exp(-c_1N^{1-\alpha}n^{\alpha})$ such that for every $z \in {\cal T}$,
\begin{equation} \label{eq:individual-net}
\bigl|\{ i : |\inr{z,X_i}| \geq c^\prime\}| \geq 2N^{1-\alpha} n^{\alpha}.
\end{equation}

The natural choice of a set ${\cal T}$ is a minimal  $\rho$-cover of $\partial L_r^\circ$ with respect to the $\ell_2$ norm. Note that $L_r^\circ ={\rm absconv}(B_1^n \cup r^{-1/2}B_2^n) \subset B_2^n$, and so there is a $\rho$-cover of the allowed cardinality for
$$
\rho \leq 5\exp(-c_2 (N/n)^{1-\alpha}),
$$
where $c_2$ is an absolute constant.

Clearly, $\{z-\pi z : z \in \partial L_r^\circ \} \subset \rho B_2^n$, and to complete the proof of Theorem \ref{thm:main} it suffices to show that with probability at least $1-2\exp(-c_3N^{1-\alpha}n^{\alpha})$
\begin{equation} \label{eq:osc}
Q=\sup_{u \in \rho B_2^n} \bigl|\{ i : |\inr{u,X_i}| \geq c^\prime/2 \}\bigr| \leq N^{1-\alpha}n^{\alpha}.
\end{equation}

To prove \eqref{eq:osc}, observe that $Q$ is the supremum of an empirical process indexed by a class of binary valued functions
$$
F = \bigl\{f_z=\IND_{\{ |\inr{z,\cdot}| \geq c^\prime/2\}} : z \in \rho B_2^n \bigr\};
$$
in particular, for every $f_z \in F$,
\begin{equation*}
\|f_z\|_{L_2} =  Pr^{1/2}(|\inr{z,X}| \geq c^\prime/2) \leq \frac{2\|\inr{z,X}\|_{L_2}}{c^\prime} \leq \frac{2\rho}{c^\prime}
= c_4(\kappa,\delta) \exp(-c_2 (N/n)^{1-\alpha}).
\end{equation*}
By Talagrand's concentration inequality for bounded empirical processes (\cite{Tal}, see also \cite{BLM}), with probability at least $1-2\exp(-t)$,
\begin{align*}
Q \lesssim & \E Q + \sqrt{t} \sqrt{N} \sup_{f_z \in F} \|f_z\|_{L_2} + t \sup_{f_z \in F} \|f_z\|_{L_\infty}
\\
\lesssim & \E Q + \sqrt{t} \sqrt{N} c_4(\kappa,\delta) \exp(-c_2 (N/n)^{1-\alpha}) + t 
\\
= & (1)+(2)+(3);
\end{align*}
Let us show that for the right choice of $t$ and $N$ large enough, $Q \leq N^{1-\alpha} n^{\alpha}$.

The required estimate on $(2)$ and $(3)$ clearly holds as long as
$$
t \lesssim_{\kappa,\delta} N^{1-\alpha} n^{\alpha} \ \ {\rm and } \ \ N \gtrsim_\alpha n.
$$
As for $\E Q$,  note that point-wise
$$
\sup_{u \in \rho B_2^n} \bigl|\{ i : |\inr{u,X_i}| \geq c^\prime/2 \}\bigr| \leq \frac{2}{c^{\prime}} \sup_{u \in \rho B_2^n} \sum_{i=1}^N |\inr{u,X_i}|.
$$
Let $(\eps_i)_{i=1}^N$ be independent, symmetric, $\{-1,1\}$-valued random variables that are independent of $(X_i)_{i=1}^N$. By the Gin\'{e}-Zinn symmetrization theorem \cite{GZ} and the contraction inequality for Bernoulli processes \cite{LT},
\begin{align*}
\E Q \leq & \frac{2}{c^{\prime}}  \E \sup_{u \in \rho B_2^n} \sum_{i=1}^N |\inr{u,X_i}|
\\
\leq & \frac{4}{c^{\prime}} \E \sup_{u \in \rho B_2^n} \sum_{i=1}^N \eps_i |\inr{u,X_i}| + \frac{2N}{c^{\prime}} \sup_{u \in \rho B_2^n} \E |\inr{u,X_i}|
\\
\leq & \frac{4}{c^{\prime}} \E \sup_{u \in \rho B_2^n} \inr{\sum_{i=1}^N \eps_i X_i,u} + \frac{2N}{c^{\prime}} \rho
\\
\leq &  \frac{4\rho}{c^{\prime}} (\sqrt{Nn} + N) \lesssim_{\kappa,\delta} N \exp(-c_2 (N/n)^{1-\alpha}),
\end{align*}
which is sufficiently small as long as $N \gtrsim_{\alpha,\kappa,\delta} n$.
\endproof

\section{Concluding Remarks}
This proof of Theorem \ref{thm:main} is based on the small-ball method and follows an almost identical path to previous results that use the method: first, one obtains an individual estimate that implies that for each $v$ in a fine-enough net, many of the values $(|\inr{X_i,v}|)_{i=1}^N$ are in the `right range'; and then, that the `oscillation vector' $(|\inr{X_i,z-v}|)_{i=1}^N$ does not spoil too many coordinates when $v$ is `close enough' to $z$. Thus, with high probability and uniformly in $z$, many of the values $(|\inr{X_i,z}|)_{i=1}^N$ are in the right range.

Having said that, there is one substantial difference between this proof and other instances in which the small-ball method had been used. Perviously, individual estimates had been obtained in the small-ball regime; here the necessary regime is different: one requires a lower estimate on the tails of marginals of $X=(\xi_i)_{i=1}^n$. And indeed, the core of the proof is the individual estimate from Theorem \ref{thm:single}, where one shows that if $\xi$ satisfies a small-ball condition and $X$ has iid coordinates distributed as $\xi$ then its marginals exhibit a `super-gaussian' behaviour at the right level.

\bibliographystyle{plain}
\bibliography{SB}

\end{document}